\documentclass[twocolumn,nameyear,amsthm]{autart}    

\usepackage{graphicx}          
\usepackage{amsmath}
\usepackage{amsfonts}

\usepackage{color} 
\usepackage{natbib}

\usepackage[colorlinks,citecolor=blue]{hyperref}

\begin{document}

\begin{frontmatter}

\title{Continuous-Time Inverse Quadratic Optimal Control Problem  } 


\author[kth,hit]{Yibei Li}\ead{yibei@kth.se},    
\author[hit]{Yu Yao}\ead{yaoyu@hit.edu.cn},               
\author[kth]{Xiaoming Hu}\ead{hu@kth.se}  

\address[kth]{Department of Mathematics, KTH Royal Institute of Technology, SE-100 44, Stockholm, Sweden}  
\address[hit]{Control and Simulation Center, Harbin Institute of Technology, 150001 Harbin, P.~R.~China}

\begin{keyword}                           
Inverse optimization; Linear quadratic problem; Linear matrix inequality.               
\end{keyword}                             

\begin{abstract}                          
In this paper, the problem of finite horizon inverse optimal control (IOC) is investigated, where the quadratic cost function of a dynamic process is required to be recovered based on the observation of optimal control sequences.
We propose the first complete result of the necessary and sufficient condition for the existence of corresponding LQ cost functions. Under feasible cases, the analytic expression of the whole solution space is derived and the equivalence of weighting matrices in LQ problems is discussed. For infeasible problems, an infinite dimensional convex problem is formulated to obtain a best-fit approximate solution with minimal control residual. And the optimality condition is solved under a static quadratic programming framework to facilitate the computation. Finally, numerical simulations are used to demonstrate the effectiveness and feasibility of the proposed methods.
\end{abstract}

\end{frontmatter}

\section{Introduction}
In recent years, the problem of inverse optimization has regained increasing popularity in the fields of robotics, economics, and bionics~\citep{Mombaur2010, Finn2016, Berret2011, Berret2016}. It has numerous varieties in different domains, such as the inverse reinforcement learning problem in machine learning \citep{hadfield2016cooperative}, and the mechanism design problem in game theory \citep{pavan2014dynamic}. In this paper we mainly focus on the problem of inverse optimal control, which is aimed at recovering the cost function of a dynamic process based on the observation of optimal actions.

The optimality principle has been investigated as an important tool to analyze natural phenomena, such as Fermat's law in optics and Lagrange dynamics in mechanics \citep{pauwels2016linear}. In the field of biology, it is also a general hypothesis that the behavior of living systems are generated based on some optimal criteria, which leads to a promising topic of inverse optimal control.
The basic question is that given a dynamic system, when we observe the optimal policy of a specific task, how can we recover the optimization criterion based on which the optimal policy is generated?
Such estimation could then help us develop a better understanding of the physical system and reproduce a similar optimal controller in other applications.
For example, inverse optimal control is a promising tool to investigate the mechanisms underlying the human locomotion and to implement them in humanoid robots \citep{mainprice2016goal}.

The problem of reconstructing cost functions has been investigated intensively. Among the existing literatures, one well-studied direction is to treat it as a parameter identification problem, where numerous numerical results have been developed. Under this situation the cost function is usually assumed to be a linear combination of certain basic functions, with the weights remaining to be identified. On one hand, in some papers like \cite{Mombaur2010} and \cite{Berret2011}, the problem is solved in a bilevel hierarchical framework and learning methods are utilized. But a forward optimal control problem has to be solved repeatedly in each inner loop to test optimality of a candidate cost function, which would lead to a computational bottleneck.
On the other hand, in \cite{hatz2012estimating}, \cite{keshavarz2011imputing}, \cite{johnson2013inverse}, \cite{Pauwels2014inverse} and \cite{pauwels2016linear}, the problem structure is better exploited and the optimal control model is characterized by its optimality conditions. Then the problem is reformulated as a residual optimization problem, where the inner loop forward optimal control problem is replaced by a set of constraints based on Karush-Kuhn-Tucker conditions or Hamilton-Jacobi-Bellman equations.

Among the various forms of the cost function, one important direction falls under the field of deterministic linear quadratic problems, which are not only well-defined but also popular for practical purposes. Some analytic results have also been obtained due to its special form.
The inverse LQ problem is first proposed by \cite{Kalman1964} for the following Linear Quadratic (LQ) optimal control problem:
\begin{equation}\label{eIntro}
\begin{aligned}
\mathop {\min }\limits_u {\kern 1pt} {\kern 1pt} {\kern 1pt} {\kern 1pt} & \int_0^{{\infty}} {\left( {{x^T}\left( t \right)Qx\left( t \right) + {u^T}\left( t \right)Ru\left( t \right) } \right)}dt \\
s.t.{\kern 1pt} {\kern 1pt} {\kern 1pt} {\kern 1pt} &\dot x\left( t \right) = Ax\left( t \right) + Bu\left( t \right)\\
&x\left( {{t_0}} \right) = {x_0}
\end{aligned}
\end{equation}

In general, given a stabilizable constant linear plant $\left( {A,B} \right)$, and a constant stabilizing feedback control law ${u^*}\left( t \right) = K{x^*}\left( t \right)$, the inverse optimal control problem is defined by two sub-problems:

(1) Existence: determine the necessary and sufficient conditions on matrices A, B and K, such that K is an optimal control law for some cost function in the form of Eq. (\ref{eIntro}).

(2) Solution: determine all $R$ and $Q$ in Eq. (\ref{eIntro}) corresponding to the same $K$.

For the infinite-time case, Kalman studied the single-input case (R=I) in frequency domain
with the return difference condition, which is then extended to the multi-input case by \cite{Anderson1989}. In time domain based on the study of matrix equations, \cite{jameson1973inverse} gives the necessary and sufficient condition to derive the solution of $R$ from the feedback matrix $K$. However, in that result the obtained $Q$ cannot be guaranteed to be constant and nonnegative. From then on, the results of Anderson and Jameson are extended and improved to derive different results for the existence problem, such as \cite{Fujii1984}, \cite{SUGIMOTO1987}, and \cite{Fujii1987}. Then in recent years, the tool of Linear Matrix Inequality (LMI) and optimization are used in \cite{Boyd1994} and \cite{Priess2015} to calculate the solutions of $Q$ and $R$.

However, on the other hand, the inverse LQ problem in finite time is still an open problem. To the best of our knowledge, there are only a few results related to this problem. In addition to the incomplete result of \cite{jameson1973inverse}, \cite{Nori2004} makes a step forward, showing that for any quadratic cost function, there exists a canonical form with a cross term such that it can generate the same optimal control. Then \cite{jean2018inverse} makes some extensions to investigate the uniqueness of the canonical form. But under this framework, the problem is reduced to a constrained parameter identification problem, which is however not easy to solve.

In this paper, the finite-time inverse LQ problem is investigated. Given the observation of an optimal feedback matrix, the necessary and sufficient condition is given for the existence of corresponding LQ cost functions by a LMI condition. For feasible problems, the analytic expression of the whole solution space is derived and the uniqueness of solutions are analyzed. On the other hand, for infeasible cases, a best-fit approximate solution is obtained, which minimizes the control residual. The main contribution of this paper is two-folded:

%

\begin{enumerate}
\item To the best of our knowledge, our result is the first attempt to give out a complete necessary and sufficient condition for the well-posedness of the inverse LQ problem, i.e. the existence of LQ cost functions. Unlike \cite{Nori2004} and \cite{jean2018inverse}, here we focus on the standard form without cross terms, which is more advantageous in its practical meaning. For feasible cases, the whole solution space is analyzed analytically, which also sheds new light on explaining the equivalence of weighting matrices in LQ problems.

\item In infeasible cases, approximate solutions are computed through a well-posed infinite dimensional convex problem, which is formulated to minimize the residual of optimal controllers. The optimality condition is derived by the primal-dual method in the form of a matrix boundary value problem (BVP) under the constraints of positive semi-definite cones. Instead of solving the BVP numerically, we transfer it into a static quadratic programming problem, which is more computationally efficient.
\end{enumerate}

The rest of the paper is organized as follows. In section 2, some preliminaries and notations are introduced. In section 3, the inverse LQ problem are formulated mathematically. The well-posedness and exact solutions of the inverse LQ problem is investigated in Section 4, while under infeasible cases an infinite-dimensional convex optimization problem is solved to obtain a best-fit approximate solution in Section 5. Numerical simulations are given in Section 6 and some concluding remarks are drawn in Section 7.




\section{Notations and Mathematical Preliminaries}
In this paper, we denote $\mathbb{R}^n$ as the space of $n$ dimensional column vector. $\mathbb{R}^{n \times n}$ denotes the space of $n \times n$ dimensional matrix. For any two matrices $X$ and $Y$, $X \succeq Y$ means $X-Y$ is positive semi-definite. We use $C\left[ 0,T \right]$ and $NBV\left[ 0,T \right]$ to denote the space of continuous functions and normalized bounded variations over $[0,T]$ respectively.
For some special matrix spaces, we use notations
\begin{equation*}
\begin{aligned}
\mathbb{S}^{n} &:=  \left\lbrace S \in \mathbb{R}^{n \times n},~S=S^{T} \right\rbrace , \\
\mathbb{S}_+^{n}  &:= \left\lbrace S \in \mathbb{S}^{n},~S \succeq 0 \right\rbrace, \\
\mathbb{C}_s^{n}\left[ 0,T \right]  &:= \left\lbrace C\left( t \right)=C^{T}\left( t \right),~C_{ij} \in C\left[ 0,T \right]  \right\rbrace, \\
NBV_s^{n}\left[ 0,T \right]  &:= \left\lbrace X\left( t \right)=X^{T}\left( t \right),~X_{ij} \in NBV\left[ 0,T \right]  \right\rbrace,
\end{aligned}
\end{equation*}
to denote the space of Hermitian matrices, the cone of positive semi-definite matrices, matrices of continuous functions, and the matrices of normalized bounded variations respectively.

The spaces $ \mathbb{S}^{n} $ and $ \mathbb{S}_+^{n} $ are Hilbert spaces, on which the inner product is defined as:
\begin{equation*}
\left\langle S_1, S_2 \right\rangle = tr\left( {S_1^T} {S_2} \right) = tr\left( {S_1} {S_2} \right),
\end{equation*}
where $ tr $ denotes the traces of two matrices.

Some matrix operators are also used in this paper. $X^\dagger$ denotes the Moore-Penrose inverse. $\otimes$ denotes the Kronecker product. $\lVert \cdot \rVert_F $ denotes the Frobenius norm of a matrix. We use $vec(\cdot)$, $vech(\cdot)$ and $mat(\cdot)$ to denote vectorization, half vectorization, and matricization respectively. Let $e_i$ be the $i-th$ canonical basis vector for $\mathbb{R}^n$.  The matrix $E_{ij} \in \mathbb{R}^{n \times n}$ has one in its $(i,j)-th$ position and zeroes elsewhere, i.e. $E_{ij}={e_i}{e_j}^T$. The column-wise block matrix $B_i \in \mathbb{R}^{n^2 \times n}$ consists of $n$ blocks of size $n \times n$, where only the $i-th$ block is an identity matrix $I_n$ and the others are all zeros. Then for any matrix $X \in \mathbb{R}^{n \times n}$ and vector $x \in \mathbb{R}^{n^2}$, the operators of vectorization and matricization can be expressed in the form of linear transmission as
\begin{equation}
\begin{aligned}
& vec(X) = \sum_{i=1}^n {B_i}X{e_i}, \\
& mat(x) = \sum_{i=1}^n {B_i}^Tx{e_i}^T.
\end{aligned}
\end{equation}

The duplication matrix $D$ and elimination matrix $ L $ are defined respectively  by
\begin{equation}
\begin{aligned}
& D^T = \sum_{n \geq i \geq j \geq 1} {u_{ij}}{vec(T_{ij})}^T, \\
& L= \sum_{n \geq i \geq j \geq 1} {u_{ij}}{vec(E_{ij})}^T ,
\end{aligned}
\end{equation}
where
\begin{equation*}
\begin{aligned}
& u_{ij} = vech(T_{ij}), \\
& T_{ij} = \left\{ \begin{array}{ll}
  E_{ii} & \textit{if}~i=j, \\
  E_{ij}+E_{ji} & \textit{otherwise}.
  \end{array} \right.
\end{aligned}
\end{equation*}

Then for any symmetric matrix $X \in \mathbb{S}^n$, there exists a linear transformation between its vectorization and half vectorization as
\begin{equation}
\begin{aligned}
& vec(X) =D vech(X), \\
& vech(X) =L vec(X).
\end{aligned}
\end{equation}

\section{Problem Formulations}
Considering the standard finite time LQ problem:
\small
\begin{equation}\label{e23}
\begin{aligned}
\mathop {\min }\limits_u {\kern 1pt} {\kern 1pt} {\kern 1pt} {\kern 1pt} &{x^T}\left( T \right)Fx\left( T \right) + \int_0^{{T}} {\left( {{x^T}\left( t \right)Qx\left( t \right) + {u^T}\left( t \right)u\left( t \right) } \right)}dt \\
s.t.{\kern 1pt} {\kern 1pt} {\kern 1pt} {\kern 1pt} &\dot x\left( t \right) = Ax\left( t \right) + Bu\left( t \right)\\
&x\left( {{t_0}} \right) = {x_0}
\end{aligned}
\end{equation}
\normalsize
where $ x \in \mathbb{R}^n $, $ u \in \mathbb{R}^m $, and $Q, F \in \mathbb{S}_+^{n}$.


Here we make the standard assumption on the system that $ \left( A,B \right) $ is controllable, $B$ has full column rank. For the forward problem, it is well-known that there exists a unique optimal feedback control that minimizes the quadratic cost function:
\begin{equation}\label{e12}
u\left( t \right) = K\left( t \right)x\left( t \right)
=- {B^T}P\left( t \right)x\left( t \right),
\end{equation}
where $ P $ is the positive semi-definite solution to the following matrix Differential Riccati Equation (DRE):
\begin{equation}\label{e13}
- \dot P = PA + {A^T}P - PB{B^T}P + Q, P(T)=F.
\end{equation}

Then the inverse optimal control problem is formulated as following:
\begin{prob}\label{prob1}
Given a controllable constant linear plant $\left( {A,B} \right)$, and an optimal feedback control law $ K\left( t \right) $, estimate the constant matrices $ Q $ and $ F $ in the quadratic cost function (\ref{e23}) such that it could generate the observed optimal controller.
\end{prob}


Here the inverse problem is investigated in two steps.
\begin{enumerate}
	\item \textit{Existence:} determine whether there exists a quadratic cost function that could generate the observed optimal controller, and whether the solution is unique.
	\item \textit{Reconstruction:} compute a best cost function under some optimal criterion if the existence problem is feasible; otherwise give an approximate solution.
\end{enumerate}

Firstly for the existence problem, \cite{jameson1973inverse} gives out the necessary and sufficient condition to recover a symmetric non-negative matrix $ P\left( t \right) $ from the feedback matrix $ K $.

\begin{prop}\label{thm1}
~Given a feedback matrix $ K(t) $, there exists a real symmetric solution $P\left(t\right)=P\left(t\right)^T$ satisfying $K\left(t\right)=-B^T P\left(t\right)$ if and only if $K\left(t\right)B$ is symmetric  and
\begin{equation}\label{e14}
rank\left(K\left(t\right)B\right)=rank\left(K\left(t\right)\right).
\end{equation}

Then all real symmetric $P\left(t\right)$ satisfying $B^{T}P\left(t\right)=-K\left(t\right)$ are presented by:
\begin{equation}\label{e15}
P\left(t\right)=-K^{T}\left(t\right)\left(K\left(t\right)B\right)^{\dagger}K\left(t\right)+Y\left(t\right),
\end{equation}
where $Y\left(t\right)$ are all real matrices that satisfy:
\begin{equation}\label{e16}
B^{T}Y\left(t\right)=0, Y\left(t\right)=Y^{T}\left(t\right).
\end{equation}

And $P\left(t\right)$ is nonnegative if and only if the eigenvalues (must be real) of $K\left(t\right)B$ are nonpositive, and $Y\left(t\right)=Y^{T}\left(t\right) \succeq 0$.

Then the matrix $Q$ can be computed by $P\left(t\right)$ through DRE, and $ F $ is determined by $ F = P\left(T\right) = P_0\left(T\right)+Y\left(T\right) $.
\end{prop}

However, the above conditions cannot guarantee a constant and nonnegative matrix $ Q $, which does not exactly solve Problem~\ref{prob1}.

Denote $P_0\left(t\right) = -K^{T}\left(t\right)\left(K\left(t\right)B\right)^{+}K\left(t\right)$.
Substituting Eq.(\ref{e15}) into the Ricatti equation, we could simplify the nonlinear constraint of Riccati equation into the following linear one, which is in fact a Lyapunov differential equation:
\begin{equation}\label{e18}
\dot Y\left(t\right) = -{A^T}Y\left(t\right) - Y\left(t\right)A - Q - G\left(t\right),
\end{equation}
where $G\left(t\right)= \dot{P_0} + {A^T}{P_0} + {P_0}A - {P_0}B{B^{T}}{P_0}$.

On the other hand, in order to get rid of the time-variant constraint of $Y(t) \succeq 0$, we notice that for the differential Ricatti equation in Eq. (\ref{e13}), given any $Q \succeq 0$, the boundary condition $P(T) \succeq 0$ is enough to guarantee the positive semi-definiteness of $P(t)$ throughout the time interval $[0,T]$. Thus in our paper the constraint $P\left(t\right) \succeq 0$ is characterized by $P(T) = Y\left( T \right)+P_0(T) \succeq 0$. Then the inverse LQ problem is reformulated as:
\begin{prop}\label{thm2}
The observed feedback control matrix $ K\left( t \right) $ is optimal to some quadratic cost function in the form of Eq.(\ref{e23}) if and only if
$K\left(t\right)B$ is symmetric with nonpositive eigenvalues and
\begin{equation}\label{e20}
rank\left(K\left(t\right)B\right)=rank\left(K\left(t\right)\right),
\end{equation}
and there exists $Q \in \mathbb{S}_+^{n}$, $ Y\left( t \right) \in \mathbb{C}_s^n\left[ {0,T} \right] $, such that:
\begin{equation}\label{e21}
\begin{aligned}
&\left\{ \begin{array}{ll}
  \dot Y\left( t \right) =  - {A^T}Y\left( t \right) - Y\left( t \right)A - Q - G\left( t \right) \\
  {B^T}Y\left( t \right) = 0 \quad \forall t \in \left[ {0,T} \right]
  \end{array} \right. .
\end{aligned}
\end{equation}
with the boundary constraint $F = Y\left( T \right)+P_0(T) \succeq 0$.
\end{prop}

Therefore, in this paper, we mainly focus on the following problem:
\begin{prob}\label{prob2}
Find $Q \in \mathbb{S}_+^{n}$, such that:
\begin{equation}\label{e22}
\begin{aligned}
\exists  Y& \left( t \right)  \in \mathbb{C}_s^n\left[ {0,T} \right],~s.t. \\
&\left\{ \begin{array}{ll}
  \dot Y\left( t \right) =  - {A^T}Y\left( t \right) - Y\left( t \right)A - Q - G\left( t \right) \\
  {B^T}Y\left( t \right) = 0 \quad \forall t \in \left[ {0,T} \right] \\
  Y\left( T \right)+P_0(T) \succeq 0
  \end{array} \right.
\end{aligned}
\end{equation}
\end{prob}

\section{Exact Solution to the Inverse Problem} \label{SecExact}
In this section the analytic solutions to Problem \ref{prob2} is investigated. For any feasible solution $Q=Q^T \geq 0$, there exists a unique solution $Y(t)$, whose expression can be computed explicitly. Then the existence problem is equivalent to the feasibility of a LMI problem. Furthermore, for feasible problems, the structure of the solution space is analyzed and an optimal solution $Q$ can be obtained through semi-definite programming (SDP) under some optimal criterion.

\subsection{Existence Problem}
The existence problem for the inverse LQ problem is studied in this part. Given the observation of an optimal controller $K(t)$, a necessary and sufficient condition is given for the existence of a corresponding quadratic cost function.
\subsubsection{Single Input Case} \label{SI}
In order to make the expressions clear and straightforward, in this part we first start with the single-input case where $m=1$. The results will be naturally extended to the multiple-input case in Section \ref{MI}.

Firstly, a basic lemma is given, which will be used throughout this section.
\begin{lem} \label{cor1}
If the original system $(A,B)$ is controllable, then the matrix
\begin{equation}\label{e8}
H=\begin{bmatrix}
I_n \otimes B^T \\
(I_n \otimes B^T )(-(I_n \otimes A^T + A^T \otimes I_n)) \\
\vdots \\
(I_n \otimes B^T )(-(I_n \otimes A^T + A^T \otimes I_n))^{n-1}
\end{bmatrix},
\end{equation}
has full column rank.
\end{lem}

\begin{pf}
Denote
\begin{equation} \label{e_barAB}
\begin{aligned}
& \tilde{A}=-{I_n} \otimes A - A \otimes {I_n}, \\
& \tilde{B}={I_n} \otimes B.
\end{aligned}
\end{equation}

We prove $ H $ has full column rank by showing that its kernel space is zero, i.e. $Ker(H)= \lbrace \theta \rbrace$. Suppose $x=\begin{bmatrix}
x_1^T & \dots & x_n^T
\end{bmatrix}^T \in Ker(H) $, where $x_i \in \mathbb{R}^n$. Then we have
\begin{equation}\label{Eqker1}
\tilde{B}^T(\tilde{A}^T)^k x = 0, \quad k=0,...,n-1.
\end{equation}

Firstly we show by induction that $B^T(A^T)^kx_i=0, \forall k=0,...,n-1, \forall i=1,...,n$.

When $k=0$, it is obvious that $B^T{x_i}=0$ for all $i=1,...,n$ since $(I_n \otimes B^T )x=0$.

Suppose $B^T\left(A^T\right)^j{x_i}=0$ for all $0 \leq j \leq k$ and we want to show $B^T\left(A^T\right)^{k+1}x_i=0$.

Through simple calculations, we have that:
\begin{equation*}
\tilde{A}^k = \begin{pmatrix}
k \\
0
\end{pmatrix}{I_n} \otimes A^k
+\begin{pmatrix}
k \\
1
\end{pmatrix} A \otimes A^{k-1}
+ \dots
+\begin{pmatrix}
k \\
k
\end{pmatrix}A^{k} \otimes {I_n}.
\end{equation*}

Then $\tilde{B}^T(\tilde{A}^T)^{k+1} x = 0$ can be rewritten as:
\begin{equation}\label{Eqker2}
\begin{aligned}
0 &= \tilde{B}^T(\tilde{A}^T)^{k+1} x \\
&=  \sum_{j=0}^{k+1} \begin{pmatrix}
k+1 \\
j
\end{pmatrix} [(A^T)^{j} \otimes {B^T}(A^T)^{k+1-j}]x.
\end{aligned}
\end{equation}

Since $B^T\left(A^T\right)^j{x_i}=0$ for all $0 \leq j \leq k$, it holds that
\begin{equation}\label{Eqker3}
[(A^T)^{j} \otimes {B^T}(A^T)^{k+1-j}]x = 0, \quad 1 \leq j \leq k+1
\end{equation}

Plugging Eq. (\ref{Eqker3}) into Eq. (\ref{Eqker2}), we can get
\begin{equation}
\begin{aligned}
& [I_n \otimes {B^T}\left( A^T \right) ^{k+1}]x = 0 \\
\Rightarrow & {B^T}\left( A^T \right) ^{k+1}{x_i} = 0, \quad i=1,...,n
\end{aligned}
\end{equation}

Denote ${\Gamma}_c= \begin{bmatrix}
{B}  & {A}{B} & \dots & {A}^{n-1}{B}
\end{bmatrix}$ as the controllability matrix of the system $\left( A, B \right)$. Then $B^T\left( A^T \right)^kx_i=0~\left(k=0,...,n-1\right)$ can be combined as ${\Gamma}_c^T {x_i} = 0$. Since $\left( A, B \right)$ is controllable, we know that ${\Gamma}_c^T$ has full column rank, thus ${x_i} = 0,~\forall i=1,...,n$.

Hence we have proved that $Ker(H)= \lbrace \theta \rbrace$, which means matrix $H$ must have full column rank. \qed
\end{pf}

For the single-input case, $m=1$ means that the matrix $H$ is a square matrix , which is nonsingular when $(A,B)$ is controllable.

\begin{lem} \label{unique}
For any $Q \in \mathbb{R}^{n \times n}$, there exists at most one $Y(t)$ that satisfies
\begin{equation} \label{e24}
\begin{aligned}
&\left\{ \begin{array}{ll}
  \dot Y\left( t \right) =  - {A^T}Y\left( t \right) - Y\left( t \right)A - Q - G\left( t \right) \\
  {B^T}Y\left( t \right) = 0
  \end{array} \right.
\end{aligned}.
\end{equation}

If $Q$ is feasible, the corresponding $Y(t)$ is uniquely determined by
\begin{equation} \label{e34}
Hvec(Y) + Nvec(Q) + f(vec(G)) = 0,
\end{equation}
where $ N $ and $f(vec(G))$ denote
\begin{equation*}
\begin{aligned}
N &=\begin{bmatrix}
0 & -\tilde{B} & -\tilde{A}\tilde{B} & \dots & -\tilde{A}^{n-2}\tilde{B}
\end{bmatrix}^T, \\
 f(vec(G)) &= \sum_{i=0}^{n-2} \begin{bmatrix}
0 \\
\vdots \\
0 \\
-\tilde{B}^T \\
\vdots \\
-\tilde{B}^T(\tilde{A}^T)^{n-2-i}
\end{bmatrix} vec(G^{(i)}).
\end{aligned}
\end{equation*}
\end{lem}

\begin{pf}
It is well-known that for a given boundary condition $Y(T)$, the solution to the first equation in Eq. (\ref{e24}) is uniquely given by $Y(t) = \int_t^T e^{A^T(\tau-t)}(Q+G(\tau))e^{A(\tau-t)} d{\tau} + e^{A^T(T-t)}Y(T)e^{A(T-t)}$. To show the Lemma is to show that for a given $Q$, there exists at most one $Y(T)$ $s.t.~B^TY(t) \equiv 0$. However, in order to derive the analytic expressions of $(Q, Y(T))$ that also satisfy $B^TY(t) \equiv 0$ for any $t \in [0,T]$, we investigate vectorized equations instead in the remaining part of this section. Suppose Eq. (\ref{e24}) has a pair of solution $(Q, Y(t))$, vectorization of the two equations leads to
\begin{equation} \label{e30}
\left\{ \begin{array}{ll}
  vec(\dot Y)=  {\tilde{A}^T}vec(Y) - vec(Q) - vec(G), \\
  {\tilde{B}^T}vec(Y) = 0,
  \end{array} \right.
\end{equation}
where $ \tilde{A} $ and $ \tilde{B} $ are defined in (\ref{e_barAB}),and '$(t)$' is omitted for the sake of brevity.

Taking $n-1$ derivatives of ${\tilde{B}^T}vec(Y)=0$ and plugging in the first equation in Eq. (\ref{e30}), it holds that
\begin{equation}
\begin{aligned}
0 &=
\begin{bmatrix}
{\tilde{B}^T}vec(Y) \\
{\tilde{B}^T}vec(\dot Y) \\
\vdots \\
{\tilde{B}^T}vec(Y^{(n-1)})
\end{bmatrix} \\
&=Hvec(Y) + Nvec(Q) + f(vec(G)),
\end{aligned},
\end{equation}
where the subscript $i$ denotes the $i-th$ derivative, and $H$ and $N$ are the same as that in Eq. (\ref{e34}).

Since $H$ has full column rank, we know that for any feasible solutions, $Y(t)$ is uniquely determined by $Q$. \qed
\end{pf}

We have shown that if the inverse LQ problem has solutions, then for any feasible $Q$, $Y(t)$ and corresponding $F$ is uniquely determined by an explicit expression in Eq. (\ref{e34}). Then $Y(t)$ can be regarded as a function of $Q$ and the necessary and sufficient condition for the consistency of Eq. (\ref{e24}) is then given only in terms of $Q$.

\begin{thm} \label{thm3}
There exists solution to Eq. (\ref{e24}) if and only if there exists $q_v \in \mathbb{R}^{n^2}$ satisfying
\begin{multline} \label{e33}
(H\tilde{A}^TH^{-1}N+H){q_v} \\
= -H\tilde{A}^TH^{-1}f(vec(G))+\dot{f}(vec(G))-Hvec(G) .
\end{multline}
If Eq. (\ref{e33}) holds, a solution $Q$ and the unique corresponding $Y(t)$ to Eq. (\ref{e24}) is given by
\begin{equation}\label{eSolVec}
\begin{aligned}
& Q=mat(q_v), \\
& Y(t) = -mat(H^{-1}(Nq_v + f(vec(G)))),
\end{aligned}
\end{equation}
where $N$ and $f(vec(G))$ are the same as that in Lemma \ref{unique}.
\end{thm}

\begin{pf}
We first prove the necessity.

If there exist solutions to Eq. (\ref{e24}), by Lemma \ref{unique} we know that Eq. (\ref{eSolVec}) holds, which coincides with Eq. (\ref{e34}).

Taking the derivative of Eq. (\ref{e34}) and plugging in the first equation in Eq. (\ref{e30}) to eliminate $vec(\dot Y)$, we get that
\begin{equation*}
\begin{aligned}
-\dot{f}(vec(G)) &= H\tilde{A}^Tvec(Y) - Hvec(Q) - Hvec(G) \\
& \begin{split}
=H\tilde{A}^T(-H^{-1}Nvec(Q)-H^{-1}f(vec(G)))\\
 - Hvec(Q) - Hvec(G),
\end{split}
\end{aligned}
\end{equation*}
which gives
\begin{multline*}
(H\tilde{A}^TH^{-1}N+H)vec(Q) \\
= -H\tilde{A}^TH^{-1}f(vec(G))+\dot{f}(vec(G))-Hvec(G) .
\end{multline*}

Next we prove the sufficiency. Given a $q_v=vec(Q)$ that satisfies Eq. (\ref{e33}), we compute a $vec(Y)$ from Eq. (\ref{e34}). Then we show that such $vec(Q)$ and  $vec(Y)$ are the solutions to Eq. (\ref{e30}).

Note that the first $n$ rows of Eq. (\ref{e34}) gives
\begin{equation}
 {\tilde{B}^T}vec(Y) = 0.
\end{equation}

Taking the derivative of Eq. (\ref{e34}), it gives
\begin{equation}
\dot{f}(vec(G))=-Hvec(\dot{Y}).
\end{equation}

Plug it into Eq. (\ref{e33}) and multiply $H^{-1}$ on the left on both sides, we get that
\begin{equation}
\begin{aligned}
-vec(\dot{Y}) &= H^{-1}\dot{f}(vec(G))  \\
&\begin{split}
= \tilde{A}^TH^{-1}(Nvec(Q)+f(vec(G))) \\
+vec(Q)+vec(G)
\end{split} \\
&= \tilde{A}^TH^{-1}(-Hvec(Y))+vec(Q)+vec(G) \\
&= -\tilde{A}^Tvec(Y)+vec(Q)+vec(G).
\end{aligned}
\end{equation}

Hence $vec(Q)$ and $vec(Y)$ given by (\ref{e33}) and (\ref{e34}) must satisfy (\ref{e30}). And the uniqueness of $vec(Y)$ is shown by Lemma \ref{unique}. \qed
\end{pf}

For convenience, in the remaining part we denote
\begin{equation*}
\begin{aligned}
& A_Q = (H\tilde{A}^TH^{-1}N+H), \\
& B_Q = -H\tilde{A}^TH^{-1}f(vec(G))+\dot{f}(vec(G))-Hvec(G).
\end{aligned}
\end{equation*}

Thus (\ref{e33}) can be denoted by
\begin{equation} \label{e36}
A_Q vec(Q) = B_Q.
\end{equation}

As $Q$ is required to be a constant matrix, a necessary condition for Problem \ref{prob2} to be feasible is that $B_Q$ is constant over $[0,T]$, and the above time-invariant linear equation has solutions.

Then it is obvious that equations Eq. (\ref{e36}) and Eq. (\ref{e34}) together with the symmetric and nonnegative constraints of $Q$ and $Y(T)+P_0(T)$ form the solutions to Problem \ref{prob2}.

Firstly, when we consider the vectorized equations (\ref{e30}), we should also guarantee that the matricization of the solutions are symmetric.
For the differential Lyapunov equation with $Q \in \mathbb{S}^n$, we know that since $G(t)$ is symmetric,  the symmetry of $Y(t)$ for $t \in [0,T]$ is equivalent to the boundary constraint $Y(T) \in \mathbb{S}^n$. The symmetry of matrices $Q$ and $Y(T)$ is guaranteed if we consider its half vectorizations $q=vech(Q)$ and $y_T=vech(Y(T))$. Then there exist symmetric solutions $Q$ and $Y(t)$ to Eq. (\ref{e30}) if and only if the system
\begin{equation} \label{SysQY}
\underbrace {
\begin{bmatrix}
A_QD & 0 \\
ND & HD
\end{bmatrix}}_{{A_s}}\begin{bmatrix}
q \\
y_T
\end{bmatrix} = \underbrace {\begin{bmatrix}
B_Q \\
-f(vec(G(T)))
\end{bmatrix}}_{{b_s}}
\end{equation}
is consistent, i.e. ${A_s}{A_s^\dagger}{b_s}={b_s}$.

Next the nonnegative constraints of $Q$ and $Y(T)+P_0(T)$ can also be rewritten as a set of linear constraint by
\begin{equation}
\begin{aligned}
& Q = mat(Dq) = \sum_{i=1}^n {B_i^T}Dq{e_i^T} \succeq 0 \\
& P(T) = P_0(T) + Y(T) \\
& \qquad ~ = P_0(T) +  \sum_{i=1}^n {B_i^T}D{y_T}{e_i^T} \succeq 0.
\end{aligned}
\end{equation}

With the above results, then we can propose the main theorem in this part. The feasibility of the inverse LQ problem (Problem \ref{prob2}) is transformed to a standard LMI problem as claimed in Theorem \ref{thmLMI}.

\begin{thm} \label{thmLMI}
Problem \ref{prob2} is feasible if and only if $B_Q$ is constant over $[0,T]$, $A_Q$ and $B_Q$ satisfy
\begin{equation} \label{e40}
{A_s}{A_s^\dagger}{b_s}={b_s},
\end{equation}
and the following LMI problem of $v=(v_1,...,v_r)^T \in \mathbb{R}^r$ is feasible
\begin{equation} \label{e42}
\left\{ \begin{array}{ll}
  Q_0 + \sum_{i=1}^r {v_i}{Q_i} \succeq 0, \\
  P_0(T) +  Y_{T0} + \sum_{i=1}^r {v_i}Y_{Ti} \succeq 0,
  \end{array} \right.
\end{equation}
where $\begin{bmatrix}
q^0 \\
y_T^0
\end{bmatrix}:={A_s^\dagger}{b_s}$, $\begin{bmatrix}
q^1 \\
y_T^1
\end{bmatrix},\cdots, \begin{bmatrix}
q^r \\
y_T^r
\end{bmatrix}$ span the null space of $A_s$, and $Q_0, Q_i, Y_{T0}, Y_{Ti} \in \mathbb{S}^n$ are defined by
\begin{equation*}
\begin{aligned}
& Q_0 = mat(Dq^0) = \sum_{i=1}^n {B_i^T}D{q^0}{e_i^T},  \\
& Q_i = mat(Dq^i) = \sum_{i=1}^n {B_i^T}D{q^i}{e_i^T}, \\
& Y_{T0} = mat(Dy_T^0) = \sum_{i=1}^n {B_i^T}D{y_T^0}{e_i^T},  \\
& Y_{Ti} = mat(Dy_T^i) = \sum_{i=1}^n {B_i^T}D{y_T^i}{e_i^T}.
\end{aligned}
\end{equation*}

\end{thm}

\begin{pf}
By the knowledge of linear algebra, we know that Eq. (\ref{e40}) is a sufficient and necessary condition for the linear equation (\ref{SysQY}) to have a solution. And if it holds, the solutions can be expressed in the form
\begin{equation}
\begin{bmatrix}
q \\
y_T
\end{bmatrix} = (A_QD)^\dagger B_Q + \sum_{i=1}^r {v_i}\begin{bmatrix}
q^i \\
y_T^i
\end{bmatrix},
\end{equation}
for any $v_i \in \mathbb{R}$, $i=1, \dots, r$.

Then it is obvious that Eq. (\ref{e42}) is equivalent to $Q \succeq 0$ and $P(T)=P_0(T)+Y(T) \succeq 0$ respectively.
Since $Q_0, Q_i, P_0(T), Y_{T0}, Y_{Ti} \in \mathbb{S}^n$, Eq. (\ref{e42}) is a standard LMI problem and can be solved easily with toolboxes like  \emph{Matlab cvx}. \qed
\end{pf}

\begin{rem}
If the LMI problem defined in Theorem \ref{thmLMI} is feasible, then at least one exact solution to the inverse LQR problem exists. Then an optimal $Q$ can be obtained under some criterion. If the LMI problem is infeasible, then the inverse problem has no solutions. In this case an approximate solution minimizing the control residual can be obtained in Section \ref{appro}.
\end{rem}



\subsubsection{Multiple Input Case} \label{MI}
In this part, the results for the single input case is extended to systems with multiple input, i.e. $m > 1$. Here some modifications are made to generalize the results in Section \ref{SI}, which in fact include the single input case as a special case.

Note that in Theorem \ref{unique} when we prove the uniqueness of $Y(t)$ for a given $Q$, we only use the fact that $H$ has full column rank. Hence the uniqueness $Y(t)$ also holds for the multi-input case.

However, for systems with $m > 1$, $H \in \mathbb{R}^{n^2m \times n^2}$ is not a square matrix. Thus for a given $vec(Q)$, we cannot use (\ref{e34}) to compute $vec(Y)$ directly without discussing the existence of solutions. As $H$ has full column rank, we can choose $n^2$ independent rows in $H$ to form a square matrix, denoted by $\bar{H}$. Then for the linear equations in (\ref{e34}), for the chosen rows in $H$, we also pick out the corresponding rows in $N$ and $f(vec(G))$, denoted by $\bar N$ and $\bar{f}(vec(G))$ respectively. Then Theorem \ref{thm3} can be generalized to account for the multiple input case.

\begin{thm} \label{thm4}
There exists solution to Eq. (\ref{e30}) if and only if there exists $q_v \in \mathbb{R}^{n^2}$ satisfying
\begin{multline} \label{e37}
(\bar{H}\tilde{A}^T\bar{H}^{-1}\bar{N}+\bar{H}){q_v} \\
= -\bar{H}\tilde{A}^T\bar{H}^{-1}\bar{f}(vec(G))+\dot{\bar{f}}(vec(G))-\bar{H}vec(G) .
\end{multline}
Then the unique $vec(Y)$ for a feasible $vec(Q)=q_v$ is given by
\begin{equation} \label{e38}
\bar{H}vec(Y) + \bar{N}vec(Q) + \bar{f}(vec(G)) = 0.
\end{equation}
\end{thm}

\begin{pf}
Since we assume that $B$ has full column rank, then we know that the rows of $\tilde{B}^T$ are linearly independent. Thus $\tilde{B}^T$ must be involved in $\bar{H}$. Then the first $nm$ rows in (\ref{e38}) also indicates that $\tilde{B}^T vec(Y) = 0$. And the remaining proof is similar to the proof of Theorem \ref{thm3}. which is omitted here. \qed
\end{pf}

Then with similar modifications, $vec(Q)$ can also be computed with Theorem \ref{thmLMI}, where $H$, $N$ and $f$ are replaced by $\bar H$, $\bar N$ and $\bar f$ respectively.

\begin{rem}
Note that when $m=1$, we have $\bar{H} = H$. And the results in this part are exactly the same as that in Section \ref{SI}. Thus we can conclude that the above results include the single input case as a special case, and are general for systems with any number of inputs.
\end{rem}

\subsection{Analysis of the Solution Space}    
In this part, the solution space of the inverse problem is analyzed when the feasible domain is non-empty. We will show that for a given optimal controller, the solution space of $Q$ to the inverse problem is a closed and bounded convex set, whose expression can be derived explicitly. 
The equivalence of quadratic cost functions are analyzed thoroughly, and the uniqueness of solutions is also discussed.

For any two linear quadratic cost functions defined by $(Q_1, F_1)$ and $(Q_2, F_2)$ respectively, we denote $P_i(t)~(i=1,2)$ as the corresponding solution to the (DRE) with $P_i(T)=F_i$. Then by simple computations, we know that the two cost functions could generate the same optimal control $K(t)$ if and only if
\begin{equation} \label{DQ_DRE}
\left\{ \begin{array}{ll}
  -{\Delta}{\dot P(t)} = A^T{\Delta}P(t)+{\Delta}P(t)A+{\Delta}Q, \\
  B^T {\Delta}P(t) = 0 \quad \forall t \in [0,T].
  \end{array} \right.
\end{equation}
where ${\Delta}P(t)=P_1(t)-P_2(t)$ and ${\Delta}Q=Q_1-Q_2$.

With similar techniques as in Section \ref{SI}, it holds that $vec({\Delta}P(t))=-H^{-1}Nvec(\Delta Q)$. Then it is obvious that $ {\Delta}P(t) $ is constant throughout the whole time interval, i.e. $ {\Delta}{\dot P(t)} = 0 $.

\begin{prop}
For a given optimal feedback matrix $K(t)$, denote the equivalent set of corresponding quadratic cost functions as $\mathbb{J}(K) =\lbrace (Q,F) \in \mathbb{S}_+^n \times \mathbb{S}_+^n \vert K(t)=DRE \lbrace Q,F \rbrace  \rbrace$. Then in $\mathbb{J}(K)$, the mapping between $Q$ and $F$ is bijective. It means that for any $Q (or~F)$ in $ \mathbb{J}(K)$, there exists exactly one $F (or~Q) \in \mathbb{S}_+^n$ such that $ (Q,F) \in \mathbb{J}(K) $.
\end{prop}
\begin{pf}
Assume that the inverse LQ problem has a feasible solution $(\bar{Q}, \bar{F}) \in \mathbb{J}(K)$. From Eq. (\ref{DQ_DRE}) we know that  $Q \in \mathbb{S}^n$ could generate the same $K(t)$ if and only if there exists $\Delta P \in \mathbb{S}^n$, such that $ {\Delta}Q = Q-\bar{Q} $ satisfies
\begin{equation} \label{DQ}
\left\{ \begin{array}{ll}
  A^T{\Delta}P+{\Delta}PA+{\Delta}Q=0, \\
  B^T {\Delta}P = 0,
  \end{array} \right.
\end{equation}
and the corresponding terminal penalty matrix for $Q$ is $F=\bar{F}+\Delta P$.

It is obvious that if $\Delta F=\Delta P=0$, then ${\Delta}Q=0$. Hence in $\mathbb{J}(K)$, the mapping from $Q$ to $F$ is injective.
On the other hand, if $\Delta Q=0$, by Eq. (\ref{DQ}) we know that ${(A+BD)^T}\Delta P - {\Delta P }(A+BD)=0 $ for any $D \in \mathbb{R}^{m \times n}$. Since $(A,B)$ is controllable, we can always find a $D$ such that $Re ~\lambda (A+BD) < 0$. Hence the Sylvester equation has a unique solution $ \Delta P = 0 $. Therefore, we have shown that $\Delta Q = 0 \Longleftrightarrow \Delta F = 0$, i.e. the mapping between $Q$ and $F$ is bijective. \qed
\end{pf}

Hence the solution space of the cost functions $\mathbb{J}(K)$ to the inverse LQ problem can be characterized by an equivalent set of matrix $Q$, which is denoted by
\begin{equation}
\mathcal{D}_{K+} = \{Q \in \mathbb{S}_+^n \vert (Q,F) \in \mathbb{J}(K) ~\textit{for some}~F \in \mathbb{S}_+^n \}.
\end{equation}

From the above proof we know that if there exists a feasible solution $\bar Q$ to the inverse LQ problem, then the set of $Q \in \mathbb{S}^n$ that could derive the same $K(t)$ is characterized by the affine manifold
\begin{equation}
\mathcal{D}_K = \{Q \vert Q=\bar{Q}+{\Delta}Q, {\Delta}Q \in \mathcal{S} \},
\end{equation}
where $\mathcal{S}$ denotes the linear subspace
\begin{equation*}
\mathcal{S} = \{{\Delta}Q \vert \exists \Delta P \in \mathbb{S}^n, s.t. {\Delta}Q~\textit{and}~{\Delta}P~\textit{satisfy}~(\ref{DQ}) \}.
\end{equation*}

Recall that we have assumed that $B$ has full column rank. Then the singular value decomposition of $ B^T $ is expressed as:
\begin{equation}
\begin{aligned}
B^T &= U*\begin{bmatrix}
\Sigma & 0
\end{bmatrix}
*V^{T}, \\
&=U*\begin{bmatrix}
\Sigma & 0
\end{bmatrix}
*\begin{bmatrix}
V_{1} & V_{2}
\end{bmatrix}^{T},
\end{aligned}
\end{equation}
where $ U \in \mathbb{R}^{m \times m}$ and $ V \in \mathbb{R}^{n \times n} $ are unitary matrices, $ V_{1} \in \mathbb{R}^{n \times m} $, $ V_{2} \in \mathbb{R}^{n \times \left(n-m\right) } $, and $ \Sigma \in \mathbb{R}^{m \times m}$ is a diagonal matrix with the diagonal entries equal to the singular values of $ B^T $.

Then it can be proved that $  B^T {\Delta}P = 0 $ if and only if
\begin{equation}
{\Delta}P = V_2 X V_2^T, \quad X \in \mathbb{S}^{n-m},
\end{equation}
which means that the linear subspace $ \mathcal{S} $ is uniquely determined by the system matrix $A$ and $B$ as
\begin{equation}
\mathcal{S} = \{{\Delta}Q \vert \Delta Q = A^TV_2 X V_2^T + V_2 X V_2^TA, X \in \mathbb{S}^{n-m} \}.
\end{equation}

It is obvious that the affine space $\mathcal{D}_K$ is non-empty with the dimension
\begin{equation}
 r = \dfrac{(n-m)(n-m+1)}{2}.
 \end{equation}


Since the solution $Q$ is also required to be nonnegative, the whole solution space of the inverse LQ problem is determined by
\begin{equation}
\mathcal{D}_{K+} = \mathcal{D}_K \cap \mathbb{S}_+^{n}.
\end{equation}

Then the uniqueness of the solution to the inverse problem can be analyzed as claimed in Theorem \ref{thm_sol}, whose proof is given in the appendix.
\begin{thm} \label{thm_sol}
Assume $\bar{Q} \in \mathbb{S}_+^n$ is a feasible solution to the inverse LQ problem, then the solution space of $Q$ has the following properties:
\begin{enumerate}
	\item $\mathcal{D}_{K+}$ is a closed and bounded convex set with dimension $r$. For any $Q_1, Q_2 \in \mathcal{D}_{K+}$ and $Q_1 \neq Q_2$, $\Delta Q = Q_1 - Q_2$ must be indefinite with eigenvalues on both sides of the imaginary axis.
	\item If $\bar{Q}=0$ is a feasible solution, then that is the unique solution to the inverse LQ problem.
	\item If $\bar{Q} \succ 0$, then there must exist infinite number of solutions.
	\item If $\bar{Q}$ is on the boundary of $\mathbb{S}_+^n$ (i.e. $rank(\bar{Q})<n$), uniqueness of the solution depends on the specific position of $\bar{Q}$. A sufficient condition for $\bar{Q}$ to be the unique solution can be given as claimed in Proposition \ref{prop_Quni}.
\end{enumerate}
\end{thm}

\begin{prop}\label{prop_Quni}
Suppose $\bar{Q} \in \partial~\mathbb{S}_+^n$ is a feasible solution to the inverse LQ problem. Let $ \{ \Delta Q_1, \dots, \Delta Q_r \} $ denote the basis of the subspace $ \mathcal{S} $. Denote $X^*$ as a non-zero optimal solution to
\begin{equation}\label{eqDQ_P}
\begin{aligned}
& \mathop {\min } {\kern 1pt} {\kern 1pt} {\kern 1pt} {\kern 1pt} tr(\bar{Q}X) \\
s.t.{\kern 1pt} {\kern 1pt} {\kern 1pt} {\kern 1pt} & tr(\Delta Q_i X) = 0,\quad (i=1,...,r)\\
& X \succeq 0
\end{aligned}
\end{equation}
 where $X^* = Vdiag(\lambda_1,\cdots,\lambda_l,0,\cdots,0)V^T$ for some unitary matrix $V$ and $l=rank(X^*)$. Define $\mathcal{T}$ as
 \begin{equation*}
\mathcal{T} = \{V \begin{bmatrix}
0 & 0 \\
0 & W
\end{bmatrix}V^T \vert W \in \mathbb{S}^{n-l} \}.
\end{equation*}
If $\mathcal{T} \cap \mathcal{S} = \{0\}$, then $\bar{Q} $ is the unique solution to the inverse LQ problem.
\end{prop}

\begin{exmp}
Here some examples are given to illustrate different structures of the solution space.

(1) Unique solution on $\partial~\mathbb{S}_+^n$.

For the finite-time LQ problem (\ref{e23}) on time interval $\left[ {0,1} \right]$, consider the system
\begin{equation*}
A = \begin{bmatrix}
2 & 1 \\
0 & -1
\end{bmatrix},
\quad
B = \begin{bmatrix}
0 \\
1
\end{bmatrix},
\quad
x\left(0\right) = \begin{bmatrix}
0 \\
0
\end{bmatrix},
\end{equation*}
with the cost function defined by
\begin{equation*}
Q = \begin{bmatrix}
4 & 2 \\
2 & 1
\end{bmatrix},
\quad
F = \begin{bmatrix}
1 & 0 \\
0 & 1
\end{bmatrix}.
\end{equation*}

Through checking the sufficient conditions in Proposition \ref{prop_Quni}, we can get that $Q$ is the unique solution to the inverse problem, which means that there exists no other $Q$ that would generate the same optimal controller as the given cost function.

(2) Infinite solutions on both $\partial~\mathbb{S}_+^n$ and $int(\mathbb{S}_+^n)$.

For the system in case (1), we choose another $Q$ as:
\begin{equation*}
Q = \begin{bmatrix}
0 & 0 \\
0 & 2
\end{bmatrix}.
\end{equation*}

Through solving the inverse problem, we can know that all the cost functions with $Q$ in
\begin{equation*}
\mathcal{D}_{K+} = \{ \begin{bmatrix}
0 & 0 \\
0 & 2
\end{bmatrix}
+ \alpha \begin{bmatrix}
4 & 1 \\
1 & 0
\end{bmatrix} \vert~ 0 \leq \alpha \leq 8 \}
\end{equation*}
 are equivalent to the given one in the sense that they lead to the same optimal controller.

(3) Infinite solutions only on $\partial~\mathbb{S}_+^n$.

Suppose the matrix parameters are chosen as
\begin{equation*}
A = \begin{bmatrix}
1 & -1 & 1 \\
0 & 2 & -1 \\
0 & 0 & 3
\end{bmatrix},
~
B = \begin{bmatrix}
1 & 0 \\
0 & 1 \\
0 & 1
\end{bmatrix},
~
Q = \begin{bmatrix}
0 & 0 & 0 \\
0 & 2 & 0 \\
0 & 0 & 1
\end{bmatrix}.
\end{equation*}

Then it can be computed that the solution space of the inverse problem is
\begin{equation*}
\mathcal{D}_{K+} = \{ \begin{bmatrix}
0 & 0 & 0 \\
0 & 2 & 0 \\
0 & 0 & 1
\end{bmatrix}
+ \alpha \begin{bmatrix}
0 & 0 & 0 \\
0 & 2 & -3 \\
0 & -3 & 4
\end{bmatrix} \vert~ -0.2 \leq \alpha \leq 10.2 \},
\end{equation*}
which lies on $\partial~\mathbb{S}_+^n$.
\end{exmp}

In general, when the problem is feasible, there always exist an infinite number of $Q$ leading to the same optimal controller. Therefore, we define an additional criteria to obtain an "optimal" $Q$ in some sense. Here we choose to minimize the conditional number of $Q$, which is always related to the problem of numerical stability \cite{cheney2012}.
Then the LMI in Eq. (\ref{e42}) can be reformulated as the following semidefinite programming (SDP) problem, which can be solved efficiently in polynomial-time.

\begin{prob}
(semidefinite programming)
\begin{equation}\label{sdp}
\begin{aligned}
& \mathop {\min }\limits_{v,\alpha} {\kern 1pt} {\kern 1pt} {\kern 1pt} {\kern 1pt} \alpha \\
s.t.{\kern 1pt} {\kern 1pt} {\kern 1pt} {\kern 1pt} & \alpha I_n \succeq Q_0 + \sum_{i=1}^r {v_i}{Q_i}\\
& Q_0 + \sum_{i=1}^r {v_i}{Q_i} \succeq 0 \\
& P_0(T) +  Y_{T0} + \sum_{i=1}^r {v_i}Y_{Ti} \succeq 0
\end{aligned}
\end{equation}
\end{prob}

\begin{rem}
From the results of \cite{ferrante2005parametrization,jean2018inverse}, we know that the optimal controller of the finite-horizon LQ problem can be uniquely parametrized by the solutions of the Algebraic Riccati Equation (ARE). Then it can be proved any two cost functions in (\ref{e23}) with different $Q$ lead to the same optimal controller if and only if their infinite-horizon counterparts also derive the same optimal feedback control matrix $K$. Therefore, the above analysis of equivalent LQ problems also applies to infinite horizon problems.
\end{rem}

\section{Approximate Solution for Infeasible Cases} \label{appro}

In this section we consider the cases where Problem \ref{prob2} is infeasible and the exact solution does not exist. For example, the optimal controller $K(t)$ might be obtained from noisy experimental data.
We want to find an optimal $Q$ that minimizes the residual error of the derived optimal controller. Here we denote $Y\left(T\right) = Y_T$, and the approximate problem is formulated in the following convex optimization framework.

\begin{prob}\label{prob3}
The approximate cost function is obtained by solving the following convex optimization problem on $ Y\left( t \right)  \in \mathbb{C}_s^n\left[ {0,T} \right] $, $Q \in \mathbb{S}^{n}$ and $Y_T \in S^{n}$:
\begin{equation}\label{e_appro}
\begin{aligned}
&\min_{Y\left( t \right) ,Q,Y_T} \dfrac{1}{2}\int_{0}^{T} \lVert B^{T}Y\left(t\right) \rVert_F^{2} dt \\
s.t.~&\dot Y\left(t\right) = -{A^T}Y\left(t\right) - Y\left(t\right)A - Q - G\left(t\right) \\
&Q \geq 0 \\
&Y_T + P_0(T) \geq 0
\end{aligned}
\end{equation}
\end{prob}

%

Note that for any feasible $Q$ and $Y_T$, $Y(t)$ can be uniquely determined from the Lyapunov differential equation on $[0,T]$. Then the above problem is
 always feasible and the finite optimal value can be reached. The optimal solution can be regarded as a best approximation of the inverse problem in the sense that it minimizes the control residual with the observed controller, i.e. $\int_{0}^{T} \lVert K^*\left(t\right)-K\left(t\right) \rVert_F^{2} dt$.

The basic idea of residual optimization comes from our previous paper \citep[see][]{li2018convex}. But in this paper, numerous extensions are investigated and the infinite-dimensional problem is completed solved through the transformation to a quadratic programming problem.

In order to solve the optimization problem in Eq. (\ref{e_appro}), the first constraint is rewritten as:
\begin{equation}
\mathcal{A}\left(Y,Q,Y_T \right) = 0,
\end{equation}
where $\mathcal{A}(Y,Q,Y_T) = Y(t) - Y_T + \int_{T}^{t} ( {A^T}Y(\tau) + Y(\tau)A + Q + G(\tau))  d\tau $ is an affine-linear operator from $\mathbb{C}_s^{n}\left[ 0,T \right] \times {\mathbb{S}^{n}} \times {\mathbb{S}^{n}}$ to $\mathbb{C}_s^{n}\left[ 0,T \right]$.


In order to solve this infinite-dimensional convex problem with the primal-dual method,  firstly the regularity condition has to be checked.
For the convex cost function in Eq.(\ref{e_appro}), the variables are optimized over an affine-linear equality constraint, and two convex inequality constraints where the positive ordering cones are closed with nonempty interiors.
A series of work has been conducted to give sufficient conditions for the strong duality in infinite dimensional spaces~\citep[see][]{jeyakumar1992generalizations,donato2011infinite,maugeri2010remarks}. For instance, based on the concept of strong quasi-relative interior, the general Slater's condition and the closed range of $\mathcal{A}$ are sufficient to guarantee the existence of Lagrange multipliers. Then by Proposition 5.1 and Theorem 5.1 in \cite{jeyakumar1992generalizations}, it is easy to check that the regularity condition holds for our problem.

In order to formulate the dual problem , firstly we define the Lagrange multiplier as
\begin{equation*}
\begin{aligned}
&\Lambda\left(t\right) \in NBV_S^{n}\left[0,T\right],\\
&\Delta_1 \in \mathbb{S}^{n}, \\
&\Delta_2 \in \mathbb{S}^{n}. \\
\end{aligned}
\end{equation*}

The algebraic dual of $ \mathbb{C}_s^{n}\left[ 0,T \right] $ is the matrix of bounded variations, which is denoted as:
\begin{equation*}
\left\langle C\left( t \right), \Lambda\left( t \right) \right\rangle = \int_{0}^{T} tr\left[ d\Lambda^{T}\left( t \right)C\left( t \right) \right],
\end{equation*}
where $ C\left( t \right) \in \mathbb{C}_s^{n}\left[ 0,T \right] $ and $ \Lambda\left( t \right) \in NBV_s^{n}\left[ 0,T \right] $.

The Lagrangian of the primal problem is calculated by
\begin{equation}
\begin{aligned}
& \begin{split}
L = \dfrac{1}{2}\int_{0}^{T} \lVert B^{T}Y\left(t\right) \rVert_F^{2} dt + \left\langle \mathcal{A}\left(Y,Q,Y_T \right),\Lambda(t)\right\rangle  \\
 - \left\langle Q ,\Delta_1\right\rangle  - \left\langle Y_T+P_0(T) ,\Delta_2\right\rangle
 \end{split}\\
& \begin{split}
~~= \dfrac{1}{2}\int_{0}^{T} \lVert B^{T}Y\left(t\right) \rVert_F^{2} dt + \int_{0}^{T} tr\left[ d{\Lambda(t)}^{T}\mathcal{A}\left(Y,Q,Y_T \right) \right] \\
 -tr\left(Q{\Delta_1} \right)-tr\left((Y_T+P_0(T)){\Delta_2} \right).
 \end{split}
\end{aligned}
\end{equation}

Then the Karush-Kuhn-Tucker (KKT) conditions for this convex problem can be given as
\begin{equation}\label{eKKT}
\begin{aligned}
&\delta L\left( {Y\left( t \right);h\left( t \right)} \right) = 0~for~any~h\left( t \right) \in \mathbb{C}_s^{n}\left[ 0,T \right], \\
&{{\partial L} \over {\partial Q}} = 0, {{\partial L} \over {\partial {Y_0}}} = 0, \\
&\dot Y\left(t\right) = -{A^T}Y\left(t\right) - Y\left(t\right)A - Q - G\left(t\right),\\
&Y\left(t\right) \geq 0,~\forall t \in \left[ 0,T \right], \\
&Q \geq 0,~  Y_T + P_0(T) \geq 0, \\
&tr\left(Q{\Delta_1} \right) = 0,~\Delta_1 \geq 0, \\
&tr\left((Y_T + P_0(T)){\Delta_2} \right) = 0,~\Delta_2 \geq 0,
\end{aligned}
\end{equation}
which is composed of the stationary conditions, primal feasibility, dual feasibility, and complementary slackness conditions.

As for the stationary condition, firstly the Gateaux derivative of $ L $ has to be computed. Since $\Lambda\left(t\right)$ is a bounded variations, we can assume that ${\Lambda}\left(0\right)=0$ without loss of generality.
Then the Gateaux derivative of $ L $ in the direction of $ h\left( t \right) \in \mathbb{C}_s^{n}\left[ 0,T \right] $ is calculated by:
\begin{equation}
\begin{aligned}
&\delta L\left( {Y\left( t \right);h\left( t \right)} \right) \\
=& \int_{0}^{T} tr\left[BB^{T}Y\left(t\right)h\left(t\right)\right]dt \\
& + \int_{0}^{T} tr\left\lbrace d{\Lambda}^{T}\left(t\right)\left[ h\left(t\right) + \int_{T}^{t} \left( A^{T}h\left(\tau\right)+h\left(\tau\right)A\right) d\tau \right]  \right\rbrace \\
=& \int_{0}^{T} tr\left[BB^{T}Y\left(t\right)h\left(t\right)\right]dt + \int_{0}^{T} tr \left[ d{\Lambda}^{T}(t)h\left(t\right) \right] \\
& - \int_{0}^{T} tr\left\lbrace  {\Lambda}^{T}\left(t\right)\left[ A^{T}h\left(t\right)+h\left(t\right)A \right] \right\rbrace dt.
\end{aligned}
\end{equation}

Since $h\left( t \right)$ is arbitrary, here we consider $h\left( t \right)$ that is vanishing at 0 and $T$. Then we have:
\begin{equation}
\begin{aligned}
0 &= \delta L\left( {Y\left( t \right);h\left( t \right)} \right) \\
& \begin{split}
=\int_{0}^{T} tr\left\lbrace \left[ BB^{T}Y\left(t\right)-{\Lambda}^{T}\left( t \right)A^{T}-A{\Lambda}^{T}\left( t \right)\right]h\left(t\right) \right\rbrace dt\\
 -\int_{0}^{T} tr\left\lbrace {\Lambda}^{T}\left(t\right)\dot{h}\left(t\right)  \right\rbrace  dt.
\end{split}
\end{aligned}
\end{equation}
which holds for any variation  $h\left( t \right) \in \mathbb{C}_s^{n}\left[ 0,T \right]$ with $h\left( 0 \right) = h\left( T \right) = 0$.

By the theory in the calculus of variation, we know that ${\Lambda}\left( t \right)$ must satisfy:
\begin{equation}
\dot{\Lambda}(t)
= -\dfrac{1}{2}( BB^{T}Y(t) + Y^{T}(t)BB^{T} ) +{\Lambda}^{T}(t)A^{T}+A{\Lambda}^{T}( t ).
\end{equation}

Then calculating the differential of $L$ w.r.t. $Q$ and $Y_T$ respectively, we have that:
\begin{equation}
\begin{aligned}
&{{\partial L} \over {\partial Q}} = -\int_{0}^{T}{\Lambda}\left(t\right)dt - \Delta_1 =  0, \\
&{{\partial L} \over {\partial Y_T}} = -\int_{0}^{T}d{\Lambda}\left(t\right) - \Delta_2 = -\Lambda\left( T \right) - \Delta_2 = 0.
\end{aligned}
\end{equation}

In order to give out the optimality condition in the form of differential equations, we denote $ \Omega\left(t\right) = \int_{0}^{t} {\Lambda}\left(\tau\right) d\tau $. Then the KKT condition Eq.(\ref{eKKT}) can be transformed to a boundary value problem:
\begin{equation}\label{KKT_ODE}
\begin{aligned}
&\dot Y\left(t\right) = -{A^T}Y\left(t\right) - Y\left(t\right)A - Q - G\left(t\right)\\
&\begin{split}
\dot{\Lambda}\left(t\right) = -\dfrac{1}{2}\left( BB^{T}Y\left(t\right) + Y^{T}\left(t\right)BB^{T} \right) \\
+{\Lambda}^{T}\left( t \right)A^{T}+A{\Lambda}^{T}\left( t \right)
\end{split} \\
&\dot {\Omega}\left(t\right) = {\Lambda}\left(t\right) \\
&\Lambda\left( 0 \right) = 0,~{\Omega}\left(0\right)=0 \\
&Q{\Omega}\left(T\right)  = 0,~(Y_T+P_0(T)){\Lambda}\left(T\right)  = 0 \\
&{\Omega}\left(T\right) \preceq 0,~Q \succeq 0, {\Lambda}\left(T\right) \preceq 0,~Y_T+P_0(T) \succeq 0
\end{aligned}
\end{equation}

The system of matrix differential equations in Eq.(\ref{KKT_ODE}) forms a boundary value problem (BVP) of symmetric matrices $Y$, $\Lambda$, $\Omega$ and $Q$ under the constraints of positive semi-definite cones. And the above system of matrix differential equations must be consistent since there always exist at least one optimal solution to Problem \ref{prob3}.


Generally speaking, handling of inequality constraints in the above BVP is non-trivial, which requires a-priori knowledge of the optimal solution structure and always suffers from significant numerical difficulty. Therefore, in the following part instead of solving Eq. (\ref{KKT_ODE}) numerically, we transform it into a static quadratic programming problem of initial conditions, which is well-defined and computationally tractable.

It is obvious that for this system of linear differential equations, $Y(T)$, $\Lambda(T)$ and $\Omega(T)$ are uniquely determined by each pair of $(Q, Y(0))$. Hence solving the boundary value problem in Eq.(\ref{KKT_ODE}) is equivalent to finding $(Q, Y(0))$ such that the boundary conditions at terminal time (last two lines in Eq.(\ref{KKT_ODE})) are satisfied.

\begin{prop}\label{lem1}
If $ M, N \in \mathbb{S}_+^{n} $ , then we have that:
\begin{equation}
tr\left( MN \right) \geq 0,
\end{equation}
and the equality holds if and only if $ MN=0 $
\end{prop}

Since $ Q, -\Omega(T), Y(T)+P_0(T), -{\Lambda}\left(T\right) \in \mathbb{S}_+^{n} $, it follows that $tr(-Q{\Omega}(T)) \succeq 0$ and $-tr((Y_T+P_0(T)){\Lambda}\left(T\right)) \succeq 0$. Then the boundary value problem is equivalent to
\begin{equation}\label{ApproSDP}
\begin{aligned}
\min_{Q,Y(0)}& -tr(Q{\Omega}(T))-tr(Y(T)+P_0(T)){\Lambda}\left(T\right)) \\
s.t.~&Q \succeq 0,~Y(T)+P_0(T) \succeq 0 \\
&{\Omega}\left(T\right) \preceq 0, {\Lambda}\left(T\right) \preceq 0
\end{aligned}
\end{equation}
where $Y(T)$, $\Lambda(T)$ and $\Omega(T)$ are determined by $Q$ and $ Y(0)$ through the differential equations.

\begin{rem}
The consistency of matrix differential equations in Eq. (\ref{KKT_ODE}) guarantees the existence of optimal solutions to the above problem, which might not be unique. Every optimal $Q$ with optimal value at zero can be regarded as a best approximation of the inverse LQ problem with minimal control residual.
\end{rem}

In order to derive the analytic solution of $Y(T)$, $\Lambda(T)$ and $\Omega(T)$, techniques of vectorization are utilized. Here we denote $y(t)=vech(Y(t))$, $\lambda(t)=vech(\Lambda(t))$, $\omega(t)=vech(\Omega(t))$, $q=vech(Q)$, $y_0=y(0)$, $g(t)=vech(G(t))$, $p_{0T}=vech(P_0(T))$ and choose the state vector as $z(t)=[y(t);\lambda(t);\omega(t)]$. Then the differential equations in Eq. (\ref{KKT_ODE}) can be rewritten as a linear system:
\begin{equation}
\dot{z} = \hat{A}z+\begin{bmatrix}
-q-g(t) \\
0 \\
0
\end{bmatrix}, \quad z(0)=\begin{bmatrix}
y_0 \\
0 \\
0
\end{bmatrix},
\end{equation}
where
\small
\begin{equation*}
\hat{A}=\begin{bmatrix}
-L(I \otimes A^T + A^T \otimes I)D & 0 & 0 \\
-\dfrac{1}{2}L(I \otimes (BB^T) + (BB^T) \otimes I)D & L(I \otimes A + A \otimes I)D & 0 \\
0 & I & 0
\end{bmatrix}.
\end{equation*}
\normalsize

For this time-invariant linear system, the state transition matrix is partitioned as
\begin{equation*}
e^{\hat{A}t}=\begin{bmatrix}
\Phi_{11}(t) & \Phi_{12}(t) &\Phi_{13}(t) \\
\Phi_{21}(t) & \Phi_{22}(t) &\Phi_{23}(t) \\
\Phi_{31}(t) & \Phi_{32}(t) &\Phi_{33}(t)
\end{bmatrix}.
\end{equation*}

Then $z(T)$ can be derived by:
\begin{equation}
z(T)=\begin{bmatrix}
y(T)\\
\lambda(T) \\
\omega(T)
\end{bmatrix}=\begin{bmatrix}
A_1\\
A_2 \\
A_3
\end{bmatrix} y_0 + \begin{bmatrix}
B_1\\
B_2 \\
B_3
\end{bmatrix} q +\begin{bmatrix}
C_1\\
C_2 \\
C_3
\end{bmatrix},
\end{equation}
where $A_i$, $B_i$ and $C_i$ denote
\begin{equation*}
\begin{aligned}
& A_i = \Phi_{i1}(T), \quad i=1,2,3 \\
& B_i = -\int_{0}^{T} \Phi_{i1}(T-s)ds, \quad i=1,2,3 \\
& C_i = -\int_{0}^{T} \Phi_{i1}(T-s)g(s)ds, \quad i=1,2,3
\end{aligned}
\end{equation*}

Denote the decision variable as $x_v=[q ; y_0]$.
Through some matrix computations, we can get that the cost function in Eq. (\ref{ApproSDP}) is equivalent to a quadratic function of $x_v$ denoted by
\begin{equation}
\begin{aligned}
& -tr(Q{\Omega}(T))-tr(Y(T)+P_0(T)){\Lambda}\left(T\right)) \\
= & -q^TD^TD{\omega}(T)-(y(T)+p_{0T})^TD^TD{\lambda}(T) \\
= &~x_v^TH_vx_v+f_v^Tx_v+g_v
\end{aligned}
\end{equation}
where $H_v \in \mathbb{R}^{n(n+1)\times n(n+1)}$, $f_v \in \mathbb{R}^{n(n+1)}$ and $g_v \in \mathbb{R}$ are determined by system matrices $A$, $B$ and observation $G(t)$.

Similar to the techniques in Section \ref{SI}, the constraints of positive semi-definiteness in Eq. (\ref{ApproSDP}) can also be transformed into a set of LMI constraints of $x_v$, which is denoted by $LMI(x_v) \succeq 0$. Then the optimization problem in Eq. (\ref{ApproSDP}) can be reformulated as the following quadratic programming problem with LMI constraints, which can be easily solved with the interior point method.
\begin{equation}
\begin{aligned}\label{ApproQP}
&\min_{x_v}~ x_v^TH_vx_v+f_v^Tx_v+g_v\\
& \begin{array}{r@{\quad}r@{}l@{\quad}l}
s.t.&LMI(x_v) \succeq 0\\
 &x_v \in \mathbb{R}^{n(n+1)}
\end{array}
\end{aligned}
\end{equation}

\section{Simulation Results}

In this section, numerical simulations are given to illustrate the proposed methods for solving the inverse LQ problem for both feasible and infeasible cases.

\subsection{Exact Solution of Feasible Cases}

Consider the following continuous time linear system:
\begin{equation*}
\dot{x} = Ax\left(t\right) + Bu\left(t\right),
\end{equation*}
where
$$
A = \begin{bmatrix}
1 & 0 & 1 \\
-2 & -3 & -1 \\
0 & 0 & 2
\end{bmatrix},
~
B = \begin{bmatrix}
1 & 0 \\
0 & 1 \\
0 & 1
\end{bmatrix},
\quad
x\left(0\right) = \begin{bmatrix}
1 \\
-0.5 \\
0
\end{bmatrix}.
$$

It is obvious that $ \left( A,B \right) $ is controllable. We consider the LQ problem in Eq.(\ref{e23}) on time interval $\left[ {0,1} \right]$ with the coefficients in the cost function as:
\begin{equation*}
\begin{aligned}
&Q_0 = \begin{bmatrix}
4 & -1 & -2 \\
-1 & 2 & -2 \\
2 & -2 & 3
\end{bmatrix}, \quad
&F_0 = \begin{bmatrix}
3 & -1 & 0 \\
-1 & 2 & -1 \\
0 & -1 & 1
\end{bmatrix}.
\end{aligned}
\end{equation*}

The forward LQ problem could be solved by the differential Riccati equation and the optimal feedback matrix $ K\left(t\right) $ is obtained with Eq.(\ref{e12}). Here we assume that $ K\left(t\right)$  is observed with no noise, which is then utilised to obtain precise solutions to the inverse problem.

With the observation of $K(t)$ we can solve the LMI optimization problems defined in Section \ref{SecExact} efficiently using the CVX toolbox in Matlab. The analytic solution of $ Q $ and $ F $ is then obtained as:
\begin{equation*}
\begin{aligned}
&Q =  \begin{bmatrix}
4 & -1 & 2 \\
-1 & 2 & -2 \\
2 & -2 & 3
\end{bmatrix}
+ \alpha \begin{bmatrix}
0 & -1 & 1 \\
-1 & -3 & 0 \\
1 & 0 & 3
\end{bmatrix} ,\\
&F =  \begin{bmatrix}
3 & -1 & 0 \\
-1 & 2 & -1 \\
0 & -1 & 1
\end{bmatrix}
+ \alpha \begin{bmatrix}
0 & 0 & 0 \\
0 & -0.5 & 0.5 \\
0 & 0.5 & -0.5
\end{bmatrix},
\end{aligned}
\end{equation*}
where $-0.49 \le {\alpha} \le 0.33$ is the freedom in the solution.

Among the feasible solutions, the optimal $Q$ with minimal conditional number can be obtained from the SDP problem in Eq. (\ref{sdp}) as
\begin{equation*}
Q =  \begin{bmatrix}
    4.0000 &  -0.5097 &   1.5097 \\
   -0.5097 &   3.4708 &  -2.0000\\
    1.5097 &  -2.0000 &   1.5292
\end{bmatrix}.
\end{equation*}

\subsection{Approximate Solution of Infeasible Cases}

Consider the same system as in the previous example. In order to illustrate infeasible cases, we suppose that the optimal control feedback $K(t)$ is measured with 20dB Gaussian white noise, i.e.
\[\bar K(t) = K(t) + w(t). \]

With the observation of $\bar K(t)$, the inverse problem is infeasible with no precise solution. In this case an approximate solution can be obtained from the quadratic programming problem defined in Eq. (\ref{ApproQP}), and the corresponding optimal cost function turns out to be:
\begin{equation*}
\begin{aligned}
&Q^* = \begin{bmatrix}
    3.9950 &  -0.8847  &  1.8707 \\
   -0.8847 &   2.3580  & -1.9989 \\
    1.8707 &  -1.9989  &  2.6441 \\
\end{bmatrix}, \\
&F^* = \begin{bmatrix}
    3.0015 &  -1.0024 &  -0.0023 \\
   -1.0024 &   2.0511 &  -1.0411 \\
   -0.0023 &  -1.0411 &   1.0351 \\
\end{bmatrix}.
\end{aligned}
\end{equation*}

\begin{figure}
  \centering
  \includegraphics[width=\hsize]{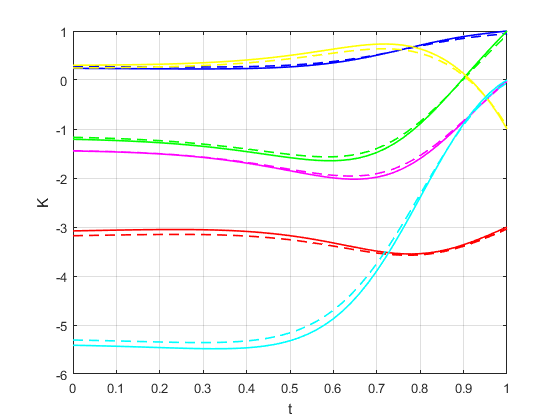}
  \caption{Curves of optimal feedback matrix $K(t)$}
  \label{fig1}
\end{figure}

\begin{figure}
  \centering
  \includegraphics[width=\hsize]{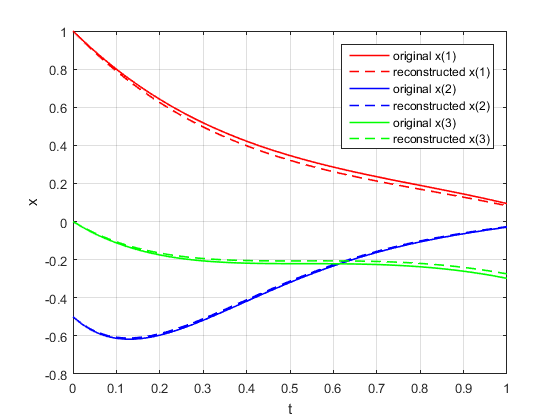}
  \caption{Curves of closed-loop state trajectory}
  \label{fig2}
\end{figure}

The curves of the optimal feedback matrix and the corresponding closed-loop state trajectory are shown in Fig.~\ref{fig1} and Fig.~\ref{fig2} respectively, where the solid line and dashed line represents the elements of original signal and reconstructed best-fit signal respectively.

The simulation shows that the recovered cost function fits the observed optimal process quite well with the optimal residual cost
\[ \int_{0}^{T} \lVert K^*\left(t\right)-K\left(t\right) \rVert_F^{2} dt = 0.0312, \]
and maximum reconstruction error
\begin{equation*}
\begin{aligned}
& \mathop {\max }\limits_t \left\| {{K^*}(t) - K\left( t \right)} \right\|_F = 0.2280, \\
& \mathop {\max }\limits_t \left\| {{x^*}(t) - x\left( t \right)} \right\| = 0.0296.
\end{aligned}
\end{equation*}

\section{Conclusions}
In this paper, we analyze the inverse LQ problem, where the existence and solutions are investigated respectively. The necessary and sufficient condition for the existence of corresponding LQ cost functions are given in the form of LMI conditions. For feasible cases, the whole solution space is shown to be a closed and bounded convex set, which is the intersection of an affine manifold and the positive semi-definite cone. And a sufficient condition for a unique cost function is also proposed. For infeasible cases, a best-fit approximate solution with minimal control residual is obtained by primal-dual method. A static quadratic programming framework is utilized to solve the optimality condition of matrix differential equations, thus improving the computational efficiency. Finally, the results of numerical simulations demonstrate the feasibility of the proposed methods and the quadratic cost function can be estimated at a high accuracy.

%


\bibliographystyle{model5-names}       
\bibliography{autosam}           

\appendix
\section{Appendix. Proofs}    
\subsection{Proof of Theorem \ref{thm_sol}}
\begin{prop}
For any $Q_1, Q_2 \in \mathcal{D}_{K+}$ and $Q_1 \neq Q_2$, $\Delta Q = Q_1 - Q_2$ must be indefinite with eigenvalues on both sides of the imaginary axis.
\end{prop}
\begin{pf}
We prove the proposition by contradiction. Suppose $\Delta Q \succeq 0$. Firstly, $B^T \Delta P = 0$ means that
\begin{equation}
ImB \subseteq Ker \Delta P.
\end{equation}

For any $x \in Ker \Delta P$, by (\ref{DQ}) we know that
\begin{equation}
x^T(A^T\Delta P+\Delta P A+\Delta Q)x = x^T{\Delta Q}x = 0.
\end{equation}

Then $\Delta Q \succeq 0$ implies ${\Delta Q}x = 0$, i.e.
\begin{equation*}
Ker \Delta P \subseteq Ker \Delta Q.
\end{equation*}

For any $x \in Ker \Delta P$, multiplying the first equation in Eq. (\ref{DQ}) by $x$ implies
\begin{equation}
 \Delta P Ax = 0~ \Rightarrow ~ Ax \in Ker \Delta P,
 \end{equation}
 which means that $Ker \Delta P$ is A-invariant.

 Thus we have that
 \begin{equation}
 \Delta P \begin{bmatrix}
 B & AB & \cdots & A^{n-1}B
 \end{bmatrix} = 0.
 \end{equation}

Since  $(A, B)$ is controllable, we must have $\Delta P = 0$, thus $\Delta Q = 0$, which is contradictory to $Q_1 \neq Q_2$.

If we assume $\Delta Q \preceq 0$, similar analysis can also show contradiction. Therefore, for any two solutions $Q_1 \neq Q_2$, $\Delta Q$ must be indefinite. \qed
\end{pf}

Consider the first property in Theorem \ref{thm_sol}. Since $\mathcal{D}_{K+}$ is obtained as the intersection of two closed convex sets, it must be convex and closed as well. And the boundedness is due to the fact that $\Delta Q$ is indefinite. Then the second and third properties can be easily derived from the first property. The forth property can be obtained by analyzing the uniqueness of solutions in LMI problems, as shown in the next brief proof.

\subsection{Proof of Proposition \ref{prop_Quni}}
Note that the dual problem of the semi-definite programming problem Eq. (\ref{eqDQ_P}) is
\begin{equation} \label{eqDQ_D}
\begin{aligned}
& \mathop {\max}\limits_{y_i} \quad 0 \\
s.t. & \sum_{i=1}^r {y_i}{\Delta Q_i} \preceq \bar Q
\end{aligned}
\end{equation}

If $\mathcal{T} \cap \mathcal{S} = \{0\}$, then $X^*$ is optimal and nondegenerate. Hence there exists a unique
optimal dual solution to (\ref{eqDQ_D}), which means that there exists a unique $\{ {{y_i}} \}_{i = 1}^r$ such that $\bar Q -\sum_{i=1}^r {y_i}{\Delta Q_i} \in \mathbb{S}_+^n$. Since $\bar Q \succeq 0$, we can conclude that $\bar Q$ is the unique solution to the inverse LQ problem.

\end{document}